\documentclass[11pt,reqno,a4paper,twoside]{amsart}
\pdfoutput=1
\usepackage{amssymb,amsmath,amscd,amstext,amsthm,amsfonts}
\usepackage[ansinew]{inputenc}
\usepackage{graphicx}
\usepackage[mathscr]{eucal}

\numberwithin{equation}{section}

\newcommand{\R}{\mathbb{R}}
\newcommand{\Pp}{\mathbb{P}}
\newcommand{\C}{\mathbb{C}}
\newcommand{\Dsk}{\mathbb{D}}
\newcommand{\N}{\mathbb{N}}

\newcommand{\Su}{\mathbb{S}}

\newcommand{\Cscr}{\mathscr{C}}
\newcommand{\Mscr}{\mathscr{M}}
\newcommand{\Iscr}{\mathscr{I}}
\newcommand{\Pscr}{\mathscr{P}}

\newcommand{\SL}{{\rm SL}}
\newcommand{\SO}{{\rm SO}}
\newcommand{\GL}{{\rm GL}}
\newcommand{\gl}{{\rm gl}}

\newcommand{\OO}{{\rm O}}
\newcommand{\diag}{{\rm diag}}
\newcommand{\Mod}[1]{\left\vert{#1}\right\vert}
\newcommand{\nrm}[1]{\left\|#1\right\|}

\newtheorem{lema}{Lemma}
\newtheorem{rmk}{Remark}
\newtheorem{prop}{Proposition}

\newtheorem{coro}{Corollary}

\newcommand{\dem}{ \par\medbreak\noindent{\bf
Proof. }\enspace}
\newcommand{\demof}[1]{ \par\medbreak\noindent{\bf
{#1}. }\enspace}
\newcommand{\cqd}{\hfill
$\sqcup\!\!\!\!\sqcap\bigskip$}

\begin{document}

\title[Formulas for Lyapunov Exponents]{Formulas for Lyapunov Exponents}

\author[A.~T.~Baraviera]{Alexandre T. Baraviera}
\address{Departamento de Matemática\\
Universidade Federal do Rio Grande do Sul\\
Porto Alegre, RS, Brasil}
\email{baravi@mat.ufrgs.br}

\author[P.~Duarte]{Pedro Duarte}
\address{Departamento de Matemática and Cmaf \\
Faculdade de Ciências\\
Universidade de Lisboa\\
Campo Grande, Edifício C6, Piso 2\\
1749-016 Lisboa, Portugal
}
\email{pduarte@ptmat.fc.ul.pt}

\date{\today}

\begin{abstract}
We derive a series summation formula  for the average logarithm norm
of the action of a matrix on the projective space.
This formula is shown to be useful to evaluate some Lyapunov exponents
 of random $\SL$-matrix cocycles, which include a special class
  for which H. Furstenberg had  provided an explicit integral formula.
\end{abstract}

\maketitle


\section{Introduction }

Lyapunov exponents  measure the rate of separation of nearby orbits in a given dynamical system.
Linear cocycles form a  class of dynamical systems where the study of Lyapunov exponents is an important and active subject.
Roughly, a linear cocycle is a dynamical system on a vector bundle
which acts linearly on fibers, over some fixed dynamics on the base of the vector bundle.
In this article we shall deal with a special class of linear cocycles where
the base dynamics is a Bernoulli shift equipped with a (constant factor) product measure,
and the fiber is the special linear group $\SL(d,\R)$.
These cocycles, that we will refer as `random linear cocycles',
were extensively studied by H. Furstenberg. See~\cite{F1},~\cite{F},~\cite{FKe},~ \cite{FKi}.
More precisely, if $G$ is some matrix Lie group, like $\SL(d,\R)$, a {\em random linear cocycle} is determined
by a probability measure $\mu$ on $G$ together with a
map $F:G^\N\times \R^d\to G^\N\times \R^d$ defined by $F(\underline{g},v)=(\sigma(\underline{g}),g_0\,v)$,
where $\underline{g}$ denotes a matrix sequence $\{ g_n\}_{n\in\N}$ and $\sigma:G^\N \to G^\N$ denotes
the shift map \, $\sigma \{ g_n\}_{n\in\N}= \{ g_{n+1}\}_{n\in\N}$.
The measure $\mu$ on $G$ determines the product Bernoulli measure $\underline{\mu}=\mu^\N$ on the bundle's base $G^\N$,
which is shift invariant. It is always assumed in this theory that $\mu$ is `integrable', which means that
$\int_G \log^+\nrm{g}\,d\mu(g)<+\infty$.
 The iterates of $F$ are given explicitly by
$F^n(\underline{g},v)=(\sigma^n\underline{g}, g_{n-1}\,\ldots\, g_1\, g_0\,v)$.
From a probabilistic point of view, this random cocycle is determined by the  independent and identically
distributed random process  $X_n:G^\N\to G$, defined by $X_n(\underline{g})=g_n=X_0(\sigma^n\underline{g})$.
Furstenberg and Kesten~\cite{FKe} have proven that the following limit, expressing
the growth rate of the product random process
$Y_n=X_{n-1} \,\ldots\, X_1 \, X_0$, always exists and is constant for $\underline{\mu}$-almost every $\underline{g}\in G^\N$,
$$ \lambda(\mu) = \lim_{n\to+\infty} \frac{1}{n}\,\log\nrm{ g_{n-1}\,\ldots\, g_1\, g_0}
=  \lim_{n\to+\infty} \frac{1}{n}\,\log\nrm{ Y_n(\underline{g})}\;. $$
The average growth rate $\lambda(\mu)$ relates to the Lyapunov exponent as follows.
Given $(\underline{g},v)\in G^\N\times\R^d$, the Lyapunov exponent
at $\underline{g}$ along $v$ is the limit
$$ \lambda_\mu(v) = \lim_{n\to+\infty} \frac{1}{n}\,\log\nrm{ g_{n-1}\,\ldots\, g_1\, g_0\,v}
=  \lim_{n\to+\infty} \frac{1}{n}\,\log\nrm{ Y_n(\underline{g})\,v} \;, $$
which exists for every $v\in\R^d-\{0\}$
and $\underline{\mu}$-almost every $\underline{g}\in G^\N$.
Moreover, this limit is constant $\underline{\mu}$-almost everywhere
and $\lambda_\mu(v)\leq \lambda(\mu)$ for every $v\in\R^d-\{0\}$, with equality for almost every $v\in\R^d$.
In fact, strict inequality can only occur for
vectors in some $\mu$-invariant  proper vector subspace $V\subset \R^d$, i.e., one which is invariant under all matrices $g$ in
the support of $\mu$. See theorem 3.5 of~\cite{FKi}.
This shows that $\lambda(\mu)$ is the largest Lyapunov exponent of the random cocycle.
Furstenberg and Kifer give a nice variational characterization of all the Lyapunov spectra in~\cite{FKi},
but we shall only deal  with the largest Lyapunov exponent $\lambda(\mu)$ here.
In~\cite{F} H. Furstenberg established the following integral formula for the Lyapunov exponent
\begin{equation}\label{Furstenberg:formula}
\lambda(\mu)  = \int_{\SL(d,\R)}  \int_{\mathbb{P}^{d-1}} \log \nrm{g\,x}\, d\nu(x)\, d\mu(g) \;,
\end{equation}
where $\Pp^{d-1}$ denotes the projective space of lines in $\R^d$, and
$\nu$ stands for any maximal $\mu$-stationary measure.
A measure $\nu$ on $\Pp^{d-1}$ is said to be {\em $\mu$-stationary} \, iff\,
$\mu\times\nu$ is an $F$-invariant measure. This amounts to say that $\nu$ is a fixed
point of the convolution  operator $P_\mu(\nu)=\mu\ast\nu=\int g_\ast\nu\,d\mu(g)$,
induced by $\mu$ on the space of Borel probability measures on $\Pp^{d-1}$.
A $\mu$-stationary measure is said to be {\em maximal} if it maximizes the left-hand-side integral
in~(\ref{Furstenberg:formula}).
See for instance section 3  of~\cite{FKi} for proofs of these facts.
Furstenberg also found very general sufficient conditions for the largest Lyapunov exponent to be strictly positive.
It is enough that the group $G$ generated by the support of $\mu$ is non compact, and no subgroup of
$G$ with finite index is reducible. A group $G$ generated by matrices in the support of $\mu$ is said to be
{\em reducible} \, iff\, there is a non trivial decomposition of $\R^d$ as a direct sum of $\mu$-invariant subspaces of $\R^d$.
See theorem 8.6 in~\cite{F}.

Furstenberg's formula~(\ref{Furstenberg:formula}) indicates a way of computing Lyapunov exponents.
But still a couple of problems persists.

\begin{enumerate}
\item  To determine the $\mu$-stationary measures  explicitly.
\item  To compute  the following integral numerically
\begin{equation}\label{RA}
R_\nu(g) = \int_{\Pp^{d-1}} \log \nrm{g\,x}\, d\nu(x)\;.
\end{equation}
\end{enumerate}

In~\cite{F} Furstenberg solves the first problem for a special class of measures on $\SL(d,\R)$.
See theorem 7.3 of~\cite{F}.
In this paper we address the second problem mainly.
First we consider the uniform Riemannian probability measure $m$ on $\Pp^{d-1}$ and prove
in section 3 that

\bigskip

\noindent
{\bf Theorem A}{\em \;
Given $g\in\SL(d,\R)$ with singular values $0 < \lambda_1\leq \lambda_2\leq \ldots\, \lambda_d$,
if $\lambda_\ast > \lambda_d/\sqrt{2}$ then the following series  converges absolutely
$$ R_m(g)=  \log \lambda_\ast -\sum_{r=1}^\infty
\frac{1}{2\,r}\,\sum_{r_1+\cdots + r_d= r} \Theta^{(d)}_{r_1,\ldots, r_d}  \,
\left(1-\frac{\lambda_1^2}{\lambda_\ast^2}\right)^{r_1}\,\ldots\,
\left(1-\frac{\lambda_d^2}{\lambda_\ast^2}\right)^{r_d} \;, $$
where
\begin{equation}\label{Theta:formula}
\Theta^{\,(d)}_{r_1,\ldots, r_{d}}  =
\frac{r!}{r_1!\cdots r_{d}!}\,\frac{(2r_1-1)!!\cdots (2r_{d}-1)!!}{d(d+2)\cdots (d+2r-2)}\;.
\end{equation}
Moreover, the coefficients $\Theta^{\,(d)}_{r_1,\ldots, r_{d}}$
form a permutation invariant  probability distribution on the finite set\,
$\Iscr_r=\{\, (r_1,\ldots, r_d)\in\N^d\,:\, r_1+\ldots + r_d=r\,\}$.
}

\bigskip

In the last section 4, we provide some applications of Theorem A.
First, we give an explicit formula for the largest Lyapunov exponent of the
random cocycles where Furstenberg was able to give explicit stationary measures.
This formula is given  in terms of the integrals $R_m(g)$,
to which we can apply Theorem A above.
In a few special cases, these formulas are used in numerical computations
of some Lyapunov exponents.
A second motivation for proving Theorem A was the role played by the integral $R_m(g)$
in the following conjecture.
For any dimension $d>2$ and every $g\in\SL(d,\R)$,
\begin{equation}\label{conj}
\int_{\SO(d)} \log \rho(k\,g)\, dk \geq R_m(g) ,
\end{equation}
where   $\rho$ stands for the spectral radius and $dk$
represents integration
with respect to the normalized Haar measure in the special orthogonal group $\SO(d)$.
This is conjectured in~\cite[Question 6.6]{BPSW}.
An analogous result is proved in~\cite{DS} for the unitary group in $\GL(d,\C)$.
Theorem A is based on the following more general result,
to be proved in section 2.

\bigskip

\noindent
{\bf Theorem B}{\em \;
Given a probability measure $\nu\in\Pscr(\Pp^{d-1})$, and
$g\in\SL(d,\R)$ with singular values $0 < \lambda_1\leq \lambda_2\leq \ldots\, \lambda_d$,
if $\lambda_\ast > \lambda_d/\sqrt{2}$ then the following series  converges absolutely
$$ R_\nu(g)=  \log \lambda_\ast -\sum_{r=1}^\infty
\frac{1}{2\,r}\,\sum_{r_1+\cdots + r_d= r} \Theta^{(d)}_{r_1,\ldots, r_d}(k,\nu)  \,
\left(1-\frac{\lambda_1^2}{\lambda_\ast^2}\right)^{r_1}\,\ldots\,
\left(1-\frac{\lambda_d^2}{\lambda_\ast^2}\right)^{r_d} \;, $$
where
$$
\Theta^{\,(d)}_{r_1,\ldots, r_{d}}(k,\nu)  =
\frac{r!}{r_1!\cdots r_{d}!}\,\int_{\Pp^{d-1}} Q_{r_1,\ldots, r_d}(k\,x)\, d\nu(x)\;,
$$
$Q_{r_1,\ldots, r_d}:\Pp^{d-1}\to\R$\, stands for the  function\,
$Q_{r_1,\ldots, r_d}(x_1,\ldots, x_d)=x_1^{2\,r_1}\, \ldots\, x_1^{2\,r_d}$, and
$k\in\OO(d,\R)$ is an orthogonal matrix such that $k g^T g k^{-1}$ is a diagonal.\\
Moreover, the coefficients $\Theta^{\,(d)}_{r_1,\ldots, r_{d}}(k,\nu)$
form a probability distribution on the finite set\,
$\Iscr_r=\{\, (r_1,\ldots, r_d)\in\N^d\,:\, r_1+\ldots + r_d=r\,\}$.
}

\bigskip

Let us remark that in theorem A the coefficients $\Theta^{\,(d)}_{r_1,\ldots, r_{d}}$:
\begin{enumerate}
\item  are given explicitly,
\item  do not depend on  the singular value decomposition of $g$, and
\item  are  invariant under permutations of the indices $(r_1,\ldots, r_d)$.
\end{enumerate}
Although more general, all these properties fail in theorem B.

\bigskip

Defining $\Theta^{\,(d)}_{r_1,\ldots, r_{d}}( \nu):=\Theta^{\,(d)}_{r_1,\ldots, r_{d}}(id,\nu)$, we have
$\Theta^{\,(d)}_{r_1,\ldots, r_{d}}(k, \nu) =\Theta^{\,(d)}_{r_1,\ldots, r_{d}}(k\,\nu)$, and the measure
$\widetilde{\nu}=k \nu$ is $\widetilde{\mu}$-stationary  with
$\widetilde{\mu}=k \mu k^ {-1}$. Hence, there is no loss of generality
in assuming that $k=id$, i.e., $g$ is a diagonal matrix, in theorem B.

\bigskip

If we have  upper bounds on
some of the coefficients $\Theta^{\,(d)}_{r_1,\ldots, r_{d}}(k \nu)$ then
we can use theorem B to obtain a corresponding lower bound for $R_\nu(g)$.
For instance, in proposition~\ref{alpha:value} we shall derive
the following universal upper bound
$$ \sum_{r_1+\cdots + r_d= r} \Theta_{r_1,\ldots, r_d} (k,\nu)\,
\left(1-\frac{\lambda_1^2}{\lambda_\ast^2}\right)^{r_1}\,\ldots\,
\left(1-\frac{\lambda_d^2}{\lambda_\ast^2}\right)^{r_d} \leq \left(1-\frac{\lambda_1^2}{\lambda_\ast^2}\right)^{r} \;. $$
Then, combining this with theorem B, we can derive the trivial lower bound
$R_\nu(g)\geq \log\lambda_1$, where $\lambda_1$ denotes the least singular value of $g$.

\bigskip

The family of functions $Q_{r_1,\ldots, r_d}:\Pp^{d-1}\to\R$, with $(r_1,\ldots, r_d)\in\N^{d}$,
separates points in $\Pp^ {d-1}$. Hence, by Stone-Weirestrass' theorem,
 the linear space spanned by these monomials is a dense subalgebra of  $\Cscr(\Pp^ {d-1})$. In particular, the measure $\nu$ is completely determined by the `momenta'
 $\Theta_{r_1,\ldots, r_d} (\nu)$. If we could devise some convergent iterative scheme to approximate
 these `$\mu$-stationary  momenta', instead of the $\mu$-stationary measure $\nu$,
 then we would apply theorem B  and get bounds on the Lyapounov exponent $R_\nu(g)$.

\bigskip

\section{A General Formula}\label{general:formula}

\bigskip

Given  integers $r_1\geq 0,\,\ldots,\, r_d\geq 0$, consider the function
$Q_{r_1,\ldots, r_d}:\Pp^{d-1}\to\R$
\begin{equation}
 Q_{r_1,\ldots, r_d}(x_1,\ldots, x_d) = x_1^{2\,r_1}\ldots x_d^{2\,r_d}\;.
\end{equation}
This is a bounded function taking values between $0$ and $1$.
Setting $r=r_1+\ldots + r_d$, the minimum value of $Q_{r_1,\ldots, r_d}$ is $0$ while the maximum value,
$(r_1/r)^{r_1}\,\ldots\, (r_d/r)^{r_d}$, is attained at the
projective points with coordinates
$\left(\pm\sqrt{r_1/r},\ldots, \pm\sqrt{r_d/r} \right)$.
Let $\Pscr(\Pp^{d-1})$ denote the space of Borel probability measures on $\Pp^{d-1}$.
Throughout this section, $\nu$ will denote any probability measure in $\Pscr(\Pp^{d-1})$.

\begin{prop}\label{int:formula}
Let $A\in\gl(d,\R)$ be a symmetric matrix of the form $A=k^{-1}\,D\,k$,
where $k^{-1}$ is an orthogonal matrix consisting of $A$'s eigenvectors and $D=\diag(\lambda_1,\ldots,\lambda_d)$
is the corresponding eigenvalue matrix. Then
\begin{equation}
 \int_{\Pp^{d-1}} \langle A\,x,x\rangle^r\, d \nu(x)
=\sum_{r_1+\cdots + r_d=r} \Theta^{\,(d)}_{r_1,\ldots, r_{d}}(k,\nu)\,
\lambda_1^{r_1}\ldots \lambda_d^{r_d}\;,
\end{equation}
where
\begin{equation}\label{Theta}
\Theta^{\,(d)}_{r_1,\ldots, r_{d}}(k,\nu)  =
\frac{r!}{r_1!\cdots r_{d}!}\,\int_{\Pp^{d-1}} Q_{r_1,\ldots, r_d}(k\,x)\, d\nu(x)\;.
\end{equation}
\end{prop}

\dem
Using the multinomial formula we get
\begin{equation*}
 \int_{\Pp^{d-1}} \left(\sum_{i=1}^d \lambda_i\,x_i^2\right)^r\, d \nu(x_1,\ldots, x_d)
=\sum_{r_1+\cdots + r_d=r} \frac{r!}{r_1!\ldots r_d!}\, \left( \int_{\Pp^{d-1}} Q_{r_1,\ldots, r_q}\,d\nu \right)  \,
\lambda_1^{r_1}\ldots \lambda_d^{r_d}\;,
\end{equation*}
which is the stated formula with $A=D=\diag(\lambda_1,\ldots,\lambda_d)$.
The general case follows in the same way but replacing $x$ by $k\,x$ before
integrating with respect to $\nu$. Note that $\langle A\,x,x\rangle = \sum_{i=1}^d \lambda_i\, (k_i\, x)^2$
with  $k\,x = (k_1\, x,\ldots, k_d\, x)$.
\cqd

\bigskip

Given a probability measure $\nu\in\Pscr(\Pp^{d-1})$,
 and integers $r_1\geq 0,\,\ldots,\, r_d\geq 0$ we shall denote the matrix function defined in~(\ref{Theta})
by \, $\Theta^\nu_{r_1,\ldots, r_{d}}:\OO(d,\R)\to\R$.

\begin{coro}\label{probability}
$\{\Theta^\nu_{r_1,\ldots, r_{d}}\}$ is a family of non-negative bounded functions such that
 $$\sum_{r_1+\cdots + r_d=r} \Theta^\nu_{r_1,\ldots, r_{d}}(k) = 1\quad
\text{ for every matrix }\; k\in\OO(d,\R)\;. $$
\end{coro}

\dem First
$$0 < \Theta^\nu_{r_1,\ldots, r_{d}} \leq \frac{r!}{r_1!\ldots r_d!}\,(r_1/r)^{r_1}\,\ldots\, (r_d/r)^{r_d}
=\frac{r!}{r^r}\, \frac{r_1^{r_1}}{r_1!}\,\ldots\, \frac{r_d^{r_d}}{r_d!}$$
because
$0 < Q_{r_1,\ldots, r_d}  \leq (r_1/r)^{r_1}\,\ldots\, (r_d/r)^{r_d}$.
By proposition~\ref{int:formula},
\begin{align*}
\sum_{r_1+\cdots + r_d=r} \Theta^\nu_{r_1,\ldots, r_{d}}(k)
&= \sum_{r_1+\cdots + r_d=r} \int_{\Pp^{d-1}}
\frac{r!}{r_1!\,\ldots\, r_d!} \,Q_{r_1,\ldots, r_{d}} (k\,x)\,d\nu(x)\\
&= \sum_{r_1+\cdots + r_d=r} \int_{\Pp^{d-1}}
\frac{r!}{r_1!\,\ldots\, r_d!} \,Q_{r_1,\ldots, r_{d}} (x)\,dk_\ast\nu(x)\\
&= \int_{\Pp^{d-1}} \langle x,x\rangle^r\, d k_\ast\nu(x)   = 1\;.
\end{align*}
\cqd

\bigskip

\demof{Proof of theorem B}
We use the following Taylor' series for the logarithm function
$$ \log x = \log x_0 -\sum_{r=1}^\infty \frac{1}{r}\,\left( 1-x_0^{-1}\,x \right)^r \;.$$
Take $x_0=\lambda_\ast^2$ \, and  \, $b=I- \lambda_\ast^{-2}\,g^T g$. \,
Note that \, $\lambda_\ast^{-2}\,\nrm{g\,x}^2=\langle \lambda_\ast^{-2}\,g^T g\,x,x \rangle$ \,
and $b=k^{-1}\,\diag\left(1-\frac{\lambda_1^2}{\lambda_\ast^2},\ldots,1-\frac{\lambda_d^2}{\lambda_\ast^2} \right)\,k$.
Hence, by proposition~\ref{int:formula},

\begin{align*}
\int_{\Pp^{d-1}} \log \nrm{g\,x}\, d\nu(x)  &= \frac{1}{2}\,\int_{\Pp^{d-1}} \log \nrm{g\,x}^2\, d\nu(x) \\
&=
\log \lambda_\ast  - \sum_{r=1}^\infty \frac{1}{2\,r}\, \int_{\Pp^{d-1}}
\left( 1- \lambda_\ast^{-2}\,\nrm{g\,x}^2 \right)^r\, d\nu(x)\\
&=
\log \lambda_\ast  - \sum_{r=1}^\infty \frac{1}{2\,r}\, \int_{\Pp^{d-1}}
 \langle b\,x,x \rangle^r\, d\nu(x)\\
&=
\log \lambda_\ast  - \sum_{r=1}^\infty \frac{1}{2\,r}\, \sum_{r_1+\cdots + r_d= r}
 \Theta^{(d)}_{r_1,\ldots, r_d} (k,\nu)\,
\left(1-\frac{\lambda_1^2}{\lambda_\ast^2}\right)^{r_1}\,\ldots\,
\left(1-\frac{\lambda_d^2}{\lambda_\ast^2}\right)^{r_d}\;.
\end{align*}
The assumption $\lambda_\ast > \lambda_d/\sqrt{2}$ implies that
$$ \alpha = \max_{1\leq i\leq d} \Mod{ 1-\frac{\lambda_i^2}{\lambda_\ast^2} }<1\;. $$
The r.h.s. converges absolutely because the absolute value series is majorated by
$$ \sum_{r=1}^\infty \frac{1}{2\,r}\, \sum_{r_1+\cdots + r_d= r} \Theta^{(d)}_{r_1,\ldots, r_d} (k,\nu)\,
\alpha^r  = \sum_{r=1}^\infty \frac{\alpha^r}{2\,r}<+\infty \;.      $$
\qed

\begin{rmk} \label{alpha}  Assuming  $\lambda_\ast > \sqrt{(\lambda_1^2 + \lambda_d^2)/2}$ we have
$$ \alpha = \max_{1\leq i\leq d} \left( 1-\frac{\lambda_i^2}{\lambda_\ast^2} \right) =1-\frac{\lambda_1^2}{\lambda_\ast^2}
\; \in (0,1)\;. $$
\end{rmk}

\bigskip

We end this section with one more remark.

\begin{prop}\label{alpha:value}
Given   $g\in\SL(d,\R)$ with singular values $0<\lambda_1\leq \ldots \leq \lambda_d$,
let $k$ be an orthogonal matrix such that
$g^T g= k^{-1}\,\diag(\lambda_1^2,\ldots,\lambda_d^2)\,k$.
If $\lambda_\ast > \sqrt{(\lambda_1^2 + \lambda_d^2)/2}$ then
$$ \sum_{r_1+\cdots + r_d= r} \Theta^{(d)}_{r_1,\ldots, r_d} (k,\nu)\,
\left(1-\frac{\lambda_1^2}{\lambda_\ast^2}\right)^{r_1}\,\ldots\,
\left(1-\frac{\lambda_d^2}{\lambda_\ast^2}\right)^{r_d} \leq \left(1-\frac{\lambda_1^2}{\lambda_\ast^2}\right)^{r} \;. $$
\end{prop}

\dem
Combine corollary~\ref{probability} with remark~\ref{alpha}.
\cqd

\bigskip
\bigskip

\section{Spherical Integrals}\label{spherical:integrals}

Theorem A follows from theorem B  using next formula.

\begin{prop} \label{monomial:integral}
Given integers $r_1\geq 0,\,\ldots,\,r_d\geq 0$,
\begin{equation}\label{Q:formula}
\int_{\Pp^{d-1}} Q_{r_1,\ldots ,r_d}\, dm  =
\frac{(2\,r_1-1)!! \ldots (2\,r_d-1)!!}{d\,(d+2)\,\ldots\, (d+2r-2)}  \;.
\end{equation}
\end{prop}

This formula  involves the concept of double factorial,
which relates with Euler's Gamma function. The double factorial is the recursive function
defined over the natural numbers by the relation $n!! = n\,(n-2)!!$
with initial conditions  $0!!=(-1)!!=1$.
The Gamma function, defined by the improper integral
$$\Gamma(x)= \int_0^\infty t^{x-1}\, e^{-t}\, dt \quad (x>0)\;,$$
is a solution of the functional equation
\begin{equation}\label{fac:eq}
\Gamma(x+1)= x\,\Gamma(x) \;.
\end{equation}
Since $\Gamma(1)= \int_0^\infty  e^{-t}\, dt=1$ \,  it follows at once that
$\Gamma(n)=(n-1)!$ for every $n\in\N$.
In other words, the Gamma function is a real analytic interpolation of the usual
factorial function over the natural numbers.
Likewise, because $\Gamma(1/2)=\sqrt{2}$ it follows easily by induction that
for every $n\in\N$,
\begin{equation}\label{double:factorial}
\Gamma\left(n+\frac{1}{2}\right) = \frac{(2n-1)!!}{2^n}\,\sqrt{\pi} \;.
\end{equation}
We refer~\cite{AAR} for a comprehensive treatment on the Gamma function.
See formula (1.1.22) there for a justification of the value $\Gamma(1/2)=\sqrt{2}$.

The Gamma function can be used to provide explicit formulas for the
volumes of spheres and balls.
Let  $\Dsk^d=\{\, x\in\R^d\,:\, \nrm{x}^2\leq 1\,\}$ be the Euclidean unit disk
and denote its volume by ${\rm V}_d$. As above, let $\Su^{d-1}$ be the Euclidean unit sphere,
i.e., the boundary of $\Dsk^d$, and denote its area by ${\rm A}_{d-1}= \int_{\Su^{d-1}} 1\,d\sigma$,
where $\sigma$ stands for the measure induced by the canonical Euclidean induced metric on $\Su^{d-1}$.
The Divergence theorem, together with a simple change of variables, may be used to
establish the following recursive relations between these volumes (see appendix A of~\cite{ABR})
\begin{equation}\label{vol:recurs}
{\rm V}_d = \frac{2\,\pi}{d}\, {\rm V}_{d-2} \quad \text{ and } \quad
{\rm A}_{d-1} = d\, {\rm V}_{d} \;.
\end{equation}
From these relations we deduce explicit formulas for the volumes of balls and spheres:
\begin{equation}\label{volumes}
{\rm V}_{d}= \frac{\pi^{d/2}}{\Gamma(1+d/2)}  \quad \text{ and } \quad
{\rm A}_{d-1}= 2\, \frac{\pi^{d/2}}{\Gamma(d/2)} \;.
\end{equation}
To see this set  ${\rm U}_d = {\pi^{d/2}}/{\Gamma(1+d/2)}$.
The functional equation~(\ref{fac:eq})
implies that ${\rm U}_d$ satisfies the same recursive equation as ${\rm V}_d$, \,
$${\rm U}_d = \frac{\pi^{d/2}}{\Gamma(1+d/2)} =  \frac{\pi}{d/2}\,\frac{\pi^{(d-2)/2}}{\Gamma(1+(d-2)/2)} =
\frac{2\,\pi}{d}\, {\rm U}_{d-2}
\;.$$
But since $\Gamma(1/2)=\sqrt{2}$, it follows that ${\rm U}_1 = 2 = {\rm V}_1$.
Also ${\rm U}_0 = 1 = {\rm V}_0$.
Hence the equality ${\rm V}_d={\rm U}_d$ holds for all $d\geq 0$.
Finally, by~(\ref{vol:recurs}) we get
$$ {\rm A}_{d-1} = d\, {\rm V}_{d} =  d\, \frac{\pi^{d/2}}{\Gamma(1+d/2)}
=  2\, \frac{\pi^{d/2}}{ \Gamma(d/2)} \;. $$

\bigskip

\demof{Proof of proposition~\ref{monomial:integral}}
We are going to reduce the integrals~(\ref{Q:formula}) to the following family of integrals
introduced in~\cite{B}
\begin{equation*}
 I_d(r_1,\ldots, r_d)  = \int_{\Dsk^{d}}
x_1^{2\,r_1}\ldots x_d^{2\,r_d}\, d x_1\,\ldots\, dx_d \;.
\end{equation*}
Using the Divergence Theorem the author deduces a recurrence formula
from which he gets the following explicit formula  for every $d\geq 2$, and every
$r_1\geq 0,\,\ldots,\,r_d\geq 0$,
\begin{equation}\label{Id}
 I_d(r_1,\ldots, r_d)  = \frac{\Gamma(r_1+\frac{1}{2})\,\ldots \, \Gamma(r_d+\frac{1}{2})}{\Gamma(r_1+\cdots + r_d +1 + \frac{d}{2})} \;.
\end{equation}
Check formula (8) of~\cite{B}.
 It is also easy to see with a change of variables' argument
that for any $r$-homogeneous function $f:\R^d\to\R$ (see Corollary 1 of~\cite{B}),
\begin{equation}\label{homog:int}
\int_{\Dsk^d} f(x_1,\ldots, x_d)\, dx_1\,\ldots\, d x_d
= \frac{1}{d+r}\,\int_{\Su^{d-1}} f(x_1,\ldots, x_d)\, d\sigma(x_1,\ldots, x_d)  \;.
\end{equation}

Now, combining~(\ref{homog:int}),~(\ref{Id}),~(\ref{volumes})  and~(\ref{double:factorial}) we get
\begin{align*}
\int_{\Pp^{d-1}} Q_{r_1,\ldots, r_d}\, dm  &=
\int_{\Su^{d-1}} x_1^{2\,r_1}\ldots x_d^{2\,r_d}\, dm(x_1,\ldots, x_d)\\
&=
\frac{1}{{\rm A}_{d-1}}\, \int_{\Su^{d-1}} x_1^{2\,r_1}\ldots x_d^{2\,r_d}\, d\sigma(x_1,\ldots, x_d)\\
&=
\frac{d+2\,r}{{\rm A}_{d-1}}\, I_d(r_1,\ldots, r_d) \\
&=
\frac{d+2\,r}{ 2\, \frac{\pi^{d/2}}{ \Gamma(d/2)}}\,
\frac{\Gamma(r_1+\frac{1}{2})\,\ldots \, \Gamma(r_d+\frac{1}{2})}{\Gamma(r +1 + \frac{d}{2})} \\
&=
\frac{d+2\,r}{ 2^{r+1}}\,
\frac{(2\,r_1-1)!!\, \ldots \, (2\,r_d-1)!!}{(r+\frac{d}{2})\,(r-1+ \frac{d}{2})\,\ldots\, \frac{d}{2} } \\
&= \frac{(2\,r_1-1)!! \ldots (2\,r_d-1)!!}{d\,(d+2)\,\ldots\, (d+2r-2)}
\end{align*}
\cqd

\bigskip

\demof{Proof of theorem A}
In view of theorem B we just have to   compute:
\begin{align*}
Q^{(d)}_{r_1,\ldots,r_d}(k,m) &= \frac{r!}{r_1!\ldots r_d!}\, \int_{\Pp^{d-1}} Q_{r_1,\ldots, r_d}(k\,x)\, dm(x)\\
&= \frac{r!}{r_1!\ldots r_d!}\, \int_{\Pp^{d-1}} Q_{r_1,\ldots, r_d}(x)\, d k_\ast m(x)
= \frac{r!}{r_1!\ldots r_d!}\, \int_{\Pp^{d-1}} Q_{r_1,\ldots, r_d} \, d  m \;,
\end{align*}
the last equality because $k_\ast m=m$, for every orthogonal matrix $k\in\OO(d,\R)$.
Combining this computation with proposition~\ref{monomial:integral} we obtain
formula~(\ref{Theta:formula}) in theorem A.
The coefficients $\Theta^{(d)}_{r_1,\ldots, r_d}$ are obviously positive rational numbers,
which by corollary~\ref{probability} form a probability distribution on the set $\Iscr_r$.
An inspection to formula~(\ref{Theta:formula}) shows these coefficients are
invariant under permutations, i.e.,
$\Theta^{(d)}_{r_1,\ldots, r_d}= \Theta^{(d)}_{r_{\pi_1},\ldots, r_{\pi_d}}$,
for every permutation $\pi$ of $\{1,\ldots, d\}$.
\cqd

\bigskip

\section{Some Applications}

Denote by $\Mscr_m$ the space of probability measures $\mu$ in $\SL(d,\R)$
that have $m$ as $\mu$-stationary measure, i.e., $\mu\ast m=m$.
The class $\Mscr_m$ is closed under orthogonal averages, i.e., if $\mu\in \Mscr_m$ then
$\int_{\SO(d,\R)} k_\ast \mu\, dm(k) \in \Mscr$.

\begin{prop} \label{Mm:lambda} For any measure $\mu\in\Mscr_m$,
its Lyapunov exponent is
$$\lambda(\mu)=\int_{\SL(d,\R)} R_m(g)\,d\mu(g)\;.$$
\end{prop}

\dem
Follows from Furstenberg integral formula~(\ref{Furstenberg:formula}).
\cqd

Theorem A can then be used to approximate this Lyapunov exponent.
A class of examples  in $\Mscr_m$ are the so called orthogonally invariant measures.
A probability $\mu\in\Pscr(\SL(d,\R))$
is said to be  {\em orthogonally invariant } if
$k_\ast\mu=\mu$ for every orthogonal matrix $k\in\SO(d,\R)$.
 We list some equivalent characterizations of orthogonally invariant
 measures.

\begin{prop}\label{orth:inv:char}
Given a measure $\mu\in\Pscr(\SL(d,\R))$, the following are equivalent:
\begin{enumerate}
\item $\mu$ is orthogonally invariant,
\item $\mu\ast\delta_p=m$, \, $\forall\, p\in\Pp^{d-1}$,
\item $\mu\ast\nu=m$, \, $\forall\, \nu\in\Pscr(\Pp^{d-1})$,
\item $\mu = m_K\ast \theta$, for some measure $\theta\in\Pscr(\SL(d,\R))$,
\end{enumerate}
where $m_K$ stands for the normalized Haar measure on $K=\SO(d,\R)$.
\end{prop}

\dem The proof is straightforward.\cqd

\bigskip

Given a matrix $g\in\SL(d,\R)$, consider the measure
\begin{equation}\label{g}
\mu= m_K\ast \delta_{g}=\int_{\SO(d,\R)} \delta_{k g}\, dm_K(k)\;.
\end{equation}

\begin{prop}\label{lambda:mu}
The measure~ (\ref{g}) is orthogonally invariant, and its
Lyapunov exponent   is $\lambda(\mu)=R_m(g)$.
\end{prop}

\dem
Since $\mu$ is orthogonally invariant we have $\mu\ast m= m$,
and hence  by proposition~\ref{Mm:lambda}
$$\lambda(\mu)=\int_{\SL} R_m(g')\,d\mu(')
= \int_{\SO} R_m(k\,g)\,dm_K(k) = R_m(g)\;, $$
because all matrices $k\,g$ have the same singular values.
\cqd

\bigskip

Consider now the matrix family \,
\begin{equation}\label{gt}
g_t=\left(\begin{array}{cc} t I_d& 0 \\
0 & t^{-1} I_d \end{array}\right)\in\SL(2d,\R)\quad (t\geq 1) \;,
\end{equation}
where $I_d$ denotes the identity $d\times d$ matrix.
Next proposition refers to the following orthogonally invariant  measure\,
$\mu_t= m_K\ast \delta_{g_t}$.

\begin{prop} For every $t>1$,  the Lyapunov exponent of $\mu_t$ is
$$\lambda_{2d}(\mu_t)=\log t - \sum_{r=1}^\infty\frac{1}{2 r}\,
\frac{d (d+2)\ldots (d+2 r-2)}{(2d)(2d+2)\ldots (2d+2r-2)}\,
\left(1-\frac{1}{t^4}\right)^r \;. $$
\end{prop}

\dem
By proposition~ \ref{lambda:mu}, \,  $\lambda_{2d}(\mu_t)= R_m(g_t)$.
Notice that matrix $g_t$ has
$d$ singular values equal to $t>1$, and $d$ singular values equal to $1/t<1$.
Thus, applying theorem A with $\lambda_\ast=t$
$$ \lambda_{2d}(\mu_t)=R_m(g_t)=\log t -\sum_{r=1}^\infty\frac{1}{2r}\, \sum_{r_1+\ldots + r_d=r}
\Theta^{(2d)}_{0,\ldots, 0,r_1,\ldots, r_d}\, \left(1-\frac{1}{t^4}\right)^r  \;.$$
Notice that if $0+\ldots + 0 + r_1+\ldots + r_d=r$,
\begin{align*}
\Theta^{(2d)}_{0,\ldots, 0,r_1,\ldots, r_d} &=
\frac{r!}{r_1!\cdots r_{d}!}\,\frac{(2r_1-1)!!\cdots (2r_{k}-1)!!}{(2d)(2d+2)\cdots (2d+2r-2)}\\
&= \Theta^{(d)}_{r_1,\ldots, r_d}\,\frac{(2d)(2d+2)\cdots (2d+2r-2)}{d(d+2)\cdots (d+2r-2)}
\end{align*}
Hence, because $\sum_{r_1+\ldots + r_d=r}\Theta^{(d)}_{r_1,\ldots, r_d}=1$,
$$
\sum_{r_1+\ldots + r_d=r}
\Theta^{(2d)}_{0,\ldots, 0,r_1,\ldots, r_d}  = \frac{(2d)(2d+2)\cdots (2d+2r-2)}{d(d+2)\cdots (d+2r-2)}\;,
$$
and we get the given formula for $\lambda_{2d}(\mu_t)$.
\cqd

\begin{coro}
The function sequence $\lambda_{2d}(\mu_t)$ increases with $d$, and for every $t\geq 1$
$$\lim_{d\to+\infty} \lambda_{2d}(\mu_t) =
\log t  + \frac{1}{2}\, \log\left(\frac{ 1 + t^4}{2\,t^4}\right)\;. $$
\end{coro}

\dem
Notice that
$$\frac{d (d+2)\ldots (d+2 r-2)}{(2d)(2d+2)\ldots (2d+2r-2)}\geq \frac{1}{2^r}\;,$$
and the right hand side decreases to $2^{-r}$ as $d$ grows to $+\infty$.
The series $\sum_{r=1}^\infty\frac{1}{2^{r+1} r}\,\left(1-\frac{1}{t^4}\right)^r$
converges absolutely and uniformly to the function
$$ g(t)= \log t - \sum_{r=1}^\infty\frac{1}{2^{r+1} r}\,\left(1-\frac{1}{t^4}\right)^r \;. $$
Because this series is essentially a geometric one we can compute its sum explicitly
$$ g(t) = \log t  + \frac{1}{2}\, \log\left(\frac{ 1 + t^4}{2\,t^4}\right)\;.$$
Then, by Lebesgue monotone convergence theorem \, $\lim_{d\to\infty} \lambda_{2d}(\mu_t)=g(t)$.
\cqd

\begin{figure}[h]
 \begin{center}
 \includegraphics*[scale=0.5]{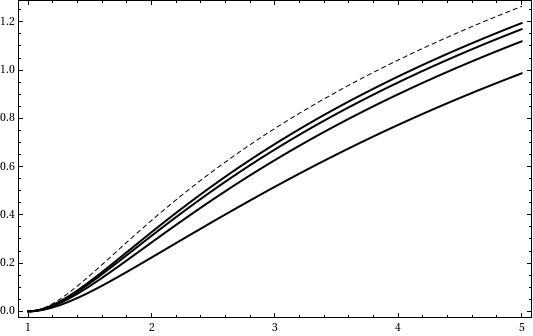}
 \end{center}
  \caption{The graphs of $t\mapsto \lambda_{2d}(\mu_t)$,
  for the half-dimensions $d=1,2,3,4$.} \label{lambda:infty}
\end{figure}

  \bigskip

This corollary shows that for large dimensions, $\lambda_{2d}(\mu_t)\approx \log t=\log \nrm{g_t}$,
which is somehow expectable since all matrices in the support of $\mu_t$ have norm $t$.

  The graphs of these functions, computed in {\em Mathematica} are depicted in figure~\ref{lambda:infty}.
  The dashed line represents the graph of $g(t)$.

\bigskip

For the following class of measures  Furstenberg was able
to give explicit stationary measures, see theorem 7.3 of \cite{F}.
Given two probability measures $\mu_1$ and $\mu_2$ in $\SL(d,\R)$, define the measure
\begin{equation}\label{m1:m2}
\mu =\mu_1\ast m_K\ast \mu_2 =  \int_{\SL}\int_{\SO}\int_{\SL}
 \delta_{g_1 k g_2}\, d\mu_1(g_1)\, dm_K(k)\, d\mu_2(g_2) \;.
\end{equation}
The $\mu$-stationary measure of $\mu$ is $\mu_1\ast m$.
In fact, by item 4. of proposition~\ref{orth:inv:char}, the measure
$m_K\ast \mu_2\ast \mu_1$ is orthogonally invariant.
Hence $(m_K\ast \mu_2\ast \mu_1)\ast m = m$  and
\begin{align*}
\mu\ast(\mu_1\ast m) &= (\mu_1\ast m_K\ast \mu_2)\ast(\mu_1\ast m)\\
&= \mu_1\ast (m_K\ast \mu_2 \ast \mu_1)\ast m \\
&= \mu_1\ast m\;,
\end{align*}
which shows that $\mu_1\ast m$ is $\mu$-stationary.

\begin{prop}\label{mu}
The Lyapunov exponent of ~ (\ref{m1:m2}) is
$$\lambda(\mu)= \int_{\SL}\int_{\SL} R_m(g_1 g_2)\, d\mu_1(g_1)\,d\mu_2(g_2) \;. $$
\end{prop}

\dem
By Furstenberg formula,
\begin{align*}
\lambda(\mu) &=\int_{\SL} \int_{\SL}  \int_{\SO}
\int_{\Pp^{d-1}} \log \nrm{g_1 k g_2 x}\, d(\mu_1\ast m)(x)\, d m_K(k)\, d\mu_1(g_1)\, d\mu_2(g_2) \\
&=\int_{\SL}\int_{\SL} \int_{\SL}  \int_{\SO}
\int_{\Pp^{d-1}} \log \nrm{g_1 k g_2 \frac{g_1' x}{\nrm{g_1' x}}}\, d m(x)\, d\mu_1(g_1')\, d m_K(k)\, d\mu_1(g_1)\, d\mu_2(g_2) \\
&=\int_{\SL} \int_{\SL} \int_{\SL} \int_{\SO}
R_m( g_1 k g_2 g_1')-R_m(g_1')\,  d m_K(k)\, d\mu_1(g_1')\, d\mu_1(g_1)\, d\mu_2(g_2)\\			
&=\int_{\SL} \int_{\SL}\int_{\SL} \int_{\SO} R_m(g_2 g_1') + R_{g_2 g_1' m}(g_1 k)  - R_m(g_1') \,  d m_K(k)
\, d\mu_1(g_1')\,  d\mu_1(g_1)\, d\mu_2(g_2) \\		
&=\int_{\SL} \int_{\SL}  R_m(g_2 g_1') - R_m(g_1')\,  d\mu_1(g_1')\, d\mu_2(g_2) \; +\\
& \qquad
\int_{\SL} \int_{\SL}\int_{\SL}\left(  \int_{\SO}  R_{g_2 g_1' m}(g_1 k)\,  d m_K(k)\right)
\, d\mu_1(g_1')\,  d\mu_1(g_1)\, d\mu_2(g_2) \\	
&=\int_{\SL} \int_{\SL}  R_m(g_2 g_1') \, d\mu_1(g_1') \, d\mu_2(g_2)  -
\int_{\SL}   R_m(g_1') \, d\mu_1(g_1')  + \int_{\SL}   R_m(g_1) \, d\mu_1(g_1)\\
&=\int_{\SL} \int_{\SL}  R_m(g_2 g_1')\,  d\mu_1(g_1')\, d\mu_2(g_2) \;.	
\end{align*}
On the fourth step we use item (b) of lemma~ \ref{Rm:lema1},
and on the sixth step we use  lemma~ \ref{Rm:lema2}.
\cqd

\bigskip

\begin{lema}\label{Rm:lema1}
Given $g',g\in\SL(d,\R)$,
\begin{enumerate}
\item[(a)] $R_m(g)\leq \log \nrm{g}$,
\item[(b)] $ R_m(g'g) = R_m(g)+ R_{g m}(g')$,
\end{enumerate}
\end{lema}

\dem
The proof of (a) is straightforward.
Item (b) holds because
\begin{align*}
 R_{g\,m}(g')  &= \int_{\Pp^{d-1}}\log\nrm{g'\,x}\, d g\,m(x)
= \int_{\Pp^{d-1}}\log\nrm{g'\,\frac{g\,x}{\nrm{g x}} }\, d m(x)\\
&= \int_{\Pp^{d-1}}\log {\nrm{g' g\,x}}\, d m(x) -
\int_{\Pp^{d-1}}\log{\nrm{g x}}\, d m(x)  \\
& =R_m(g' g) -R_m(g)\;.
\end{align*}
\cqd

\begin{lema}\label{Rm:lema2}
Given $g\in\SL(d,\R)$ and any probability measure $\nu\in\Pscr(\Pp^{d-1})$,
$$ \int_{\SO(d,\R)} R_\nu(g k ) \, d m_K(k)= R_m(g) \;.$$
\end{lema}

\dem The measure $m_K$  is orthogonally invariant.
This  because $m_K=m_K\ast \delta_I$,
where $I$ denotes the identity in $\SL(d,\R)$,
by item 4. of proposition~\ref{orth:inv:char}.
Then, by item 2. of the same proposition,
$m_K\ast \delta_x = m$. Hence
\begin{align*}
\int_{\SO} R_\nu(g k ) \, d m_K(k) &=
\int_{\SO}\int_{\Pp^ {d-1}} \log {\nrm{g k x}}\, d\nu(x) \, d m_K(k)\\
&=
\int_{\Pp^ {d-1}} \int_{\SO} \log {\nrm{g k x}} \, d m_K(k) \, d\nu(x)\\
 &=  \int_{\Pp^{d-1}} \int_{\Pp^{d-1}}\log\nrm{g z}\, d  (m_K\ast \delta_x) (z)\, dm(x)\\
 &= \int_{\Pp^{d-1}} \int_{\Pp^{d-1}}\log\nrm{g z}\, d  m(z)\, d\nu(x) \\
 & =  \int_{\Pp^{d-1}} R_{m}(g) \, d\nu(x) =  R_{m}(g)\;.
\end{align*}
\cqd

\bigskip

We consider now a special subclass of the previous.
Given two matrices $g',g\in\SL(d,\R)$  define the measure
\begin{equation}\label{g'g}
\mu_{g',g}= \delta_{g'}\ast m_K\ast \delta_g = \int_{\SO(d,\R)} \delta_{g' k g}\, dm_K(k)\; .
\end{equation}

\begin{coro}\label{mu:g',g}
The  Lyapunov exponent of the measure $\mu_{g',g}$ is
$$\lambda(\mu_{g',g})=  R_m(g' g) \;. $$
\end{coro}

For example, $\lambda(\mu_{g^{-1},g})=R_m(I)=0$ but the measure $\mu_{g^{-1},g}$
is supported on a compact group $g^{-1}\,\SO(d,\R)\,g$, and hence should have zero
Lyapunov exponent.

\bigskip

Consider now the measure $\mu_{g_s,g_t}$, where $g_s, g_t$ are matrices as defined in~ (\ref{gt}).
By corollary \ref{mu:g',g},
$\lambda(\mu_{g_s,g_t})=R_m(g_{s}\, g_{t})=R_m(g_{s\,t}) $. Hence
\begin{coro}
$\displaystyle \lim_{d\to+\infty} \lambda_{2d}(\mu_{g_s,g_t}) = \log (s\,t)
  - \frac{1}{2}\, \log\left(\frac{2\,s^4\,t^4}{ 1 + s^4 t^4}\right)$\;.
\end{coro}

Notice that  $\log(s\,t)=\max\{ \,\log\nrm{g\,x}\,:\, g\in \mbox{supp}(\mu_{g_s,g_t})\,,\; x\in\Pp^ {d-1}\,\}$,
  the norm of matrices in the support of $\mu_{g_s,g_t}$ is not constant and
 $\lambda(\mu_{g_s,g_t})$ is some kind of average of the logarithms
$\log\,\nrm{g x}$, with $ g\in \mbox{supp}(\mu_{g_s,g_t})$ and $x\in\Pp^ {d-1}$.
From this we conclude that large dimensions bring the average $\lambda(\mu_{g_s,g_t})$
closer to its maximum possible value, $\log(s\,t)$, provided  $s\,t$ is large.

A similar conclusion, assuming conjecture~(\ref{conj}) to hold, is that
for large $t>1$ and large dimension $d$,
$$ \log t=\log\nrm{g_t}\geq \int_{\SO} \log \rho(k\,g_t)\, dm_K(k) \geq R_m(g_t)\approx \log t \;. $$
Again, this shows that  large dimensions bring the average $\int_{\SO} \log \rho(k\,g_t)\, dm_K(k)$
close to the maximum value $\log\nrm{g_t}$.

\bigskip

\section*{Acknowledgements}

This work was partially supported by Funda\c c\~ao para a Ci\^encia e a Tecnologia through the project ``Randomness in Deterministic Dynamical Systems and Applications''
ref.   \\ PTDC/MAT/105448/2008.


\thispagestyle{empty}


\begin{thebibliography}{ABCDE 5}

\addcontentsline{toc}{part}{References}

\bibitem{AAR} G. Andrews, R. Askey, R. Roy,  {\em Special Functions. }
 Encycolpedia of Mathematics and its Applications, (1999) Cambridge University Press.

\bibitem{AB}  A. Avila, J. Bochi, {\em Lyapunov Exponents. }
Lecture Notes, School and Workshop on Dynamical Systems, ICTP-Trieste


\bibitem{AB2}  A. Avila, J. Bochi, {\em A formula with some applications to the theory of Lyapunov exponents. }\\
Israel J. Math., V. 131, (2002), pp 125-137

\bibitem{ABR} S. Axler, P. Bourdon, W. Ramey, {\em Harmonic Function Theory },
Graduate Texts in Mathematics 137, (2001) Springer Verlag


\bibitem{B} J. Baker {\em  Integration Over Spheres and the Divergence Theorem for Balls. }
 The American Mathematical Monthly, Vol. 104, No. 1. (Jan., 1997), pp. 36-47.

\bibitem{BPSW} K.~Burns, C.~Pugh, M.~ Shub and A.~Wilkinson, {\em Recents
results about stable ergodicity}, J. Stat. Phys. {\bf 113} (2003), pp 85-149.

\bibitem{DS} J.P.~Dedieu and M.~Shub, {\em On random and mean exponents
for unitarily invariant
probability measures on $\GL(n,\C)$},  in "Geometric Methods in Dynamical
Systems (II)-Volume in honor of Jacob Palis",
Asterisque {\bf 287} (2003), pp 1-18.


\bibitem{F1} H. Furstenberg, {\em A Poisson Formula for Semi-Simple Lie Groups. }
Annales of Mathematics, Second Series, Vol. 77, No. 2, (1963), pp 335-386


\bibitem{F} H. Furstenberg, {\em Noncommuting Random Products. }
Transactions of the American Mathematical Society, Vol. 108, No. 3, (1963), pp
377-428


\bibitem{FKe}  H. Furstenberg, H. Kesten, {\em Products of Random Matrices},
The Annals of Mathematical Statistics
Vol. 31, No. 2 (1960), pp. 457-469.


\bibitem{FKi} H. Furstenberg and Y. Kifer, {\em Random Matrix Products And Measures
On Projective Spaces}, Israel Journal Of Mathematics Vol. 46, (1983) pp. 19-20.


\end{thebibliography}
\end{document}